\documentclass[b5paper,10pt,abstractoff,DIV=calc,headings=normal,twoside]{scrartcl}
\usepackage[english]{babel}
\usepackage{lmodern}
\usepackage{microtype}
\usepackage{amsmath,amssymb,amsthm,amsfonts,graphicx,etoolbox,scrlayer-scrpage}
\usepackage[noblocks]{authblk}
\makeatletter
\patchcmd{\@maketitle}{\huge}{\Large}{}{}
\patchcmd{\abstract}{\quotation}{}{}{}
\AtBeginEnvironment{abstract}{\noindent\ignorespaces}
\AtEndEnvironment{abstract}{\par\mbox{}}
\newcommand{\shortauthor}{}
\newcommand{\shorttitle}{\@title}
\makeatother
\setkomafont{author}{\small}
\setkomafont{title}{\rmfamily\bfseries}
\setkomafont{disposition}{\rmfamily\bfseries}

\def\AMS#1{\par\noindent \textbf{AMS subject classification: }#1\par}

\newcommand{\keywords}[1]{\par\noindent\textbf{Keywords: }#1}
\rehead[]{\shortauthor}\lohead[]{\shorttitle}\lehead[]{PPP}\rohead[]{PPP}
\ofoot[]{}\ifoot[]{}\recalctypearea
\theoremstyle{plain}
\newtheorem{theorem}{Theorem}

\newtheorem{lemma}{Lemma}
\theoremstyle{definition}
\newtheorem{definition}{Definition}
\newtheorem{conjecture}{Conjecture}
\theoremstyle{remark}

\renewenvironment{abstract}{\bigskip\noindent\begin{minipage}{\textwidth}\setlength{\parindent}{15pt}\paragraph{Abstract:}}{\end{minipage}}

\begin{document}


\renewcommand{\shortauthor}{A. Basse-O'Connor and T. Overgaard and M. Skj{\o}tt}

\title{Some Results on Random Mixed {SAT} Problems}

\author[1]{Andreas Basse-O'Connor}
\author[1]{Tobias Overgaard\thanks{Corresponding author: tlo@math.au.dk}}
\author[1]{Mette Skj{\o}tt}
\affil[1]{Department of Mathematics, Aarhus University}

\maketitle

\begin{abstract}
In this short paper we present a survey of some results concerning the random SAT problems. To elaborate, the Boolean Satisfiability (SAT) Problem refers to the problem of determining whether a given set of \(m\) Boolean constraints over \(n\) variables can be simultaneously satisfied, i.e.\ all evaluate to \(1\) under some interpretation of the variables in \(\{ 0,1\}\). If we choose the \(m\) constraints i.i.d.\ uniformly at random among the set of disjunctive clauses of length \(k\), then the problem is known as the random \(k\)-SAT problem. It is conjectured that this problem undergoes a structural phase transition; taking \(m=\alpha n\) for \(\alpha>0\), it is believed that the probability of there existing a satisfying assignment tends in the large \(n\) limit to \(1\) if \(\alpha<\alpha_\mathrm{sat}(k)\), and to \(0\) if \(\alpha>\alpha_\mathrm{sat}(k)\), for some critical value \(\alpha_\mathrm{sat}(k)\) depending on \(k\). We review some of the progress made towards proving this and consider similar conjectures and results for the more general case where the clauses are chosen with varying lengths, i.e.\ for the so-called random mixed SAT problems.
\end{abstract}


\smallskip

\AMS{Primary: 60K35; secondary: 82B44, 68R99.} 


\section{Introduction}

The present paper  presents a survey of conjectures and results on the phase transitions for the random SAT problems with a focus on mixed formulas.  The Boolean Satisfiability Problem (abbreviated: the SAT problem) lies at the heart of the famous \(P\) vs.\ \(NP\) Millennium Prize Problem. Indeed, the existence of a  polynomial time algorithm for deciding satisfiability is equivalent to \(P=NP\) (see \cite{Cook1971}). In addition to the theoretical importance of the problem, the SAT problem also holds tremendous practical relevance as it arises in many applied contexts, such as in artificial intelligence, electronic design automation, bioinformatics, and more; see \cite{Silva2008}. Motivated mainly by these practical aspects, much work has since the '90s been put into finding optimized algorithms for solving the SAT problem that now perform much better than previously thought possible (cf.\ the \(P\neq NP\) conjecture). The apparent discrepancy in the difficulty of solving SAT problem instances has motivated the study of the \emph{typical} structure of the SAT problem, and we will discuss this approach further in the following. 

\section{Preliminaries}

The SAT problem asks  whether a  Boolean function \(f:\{ 0,1\}^n\to\{ 0,1\}\) on \emph{conjunctive normal form} (abbreviated CNF)  can attain the value \(1\); in the affirmative case, \(f\) is said to be \emph{satisfiable}. CNF means that \(f\) is written as a conjunction of disjunctions, or more precisely that \(f=C_1\land\dots\land C_m\), where \(\land\) denotes logical AND, and where each of the functions \(C_j:\{ 0,1\}^n\to\{ 0,1\}\) is a disjunctive \emph{clause}, i.e.\ given as the logical OR of \(k\) variables or their negation. The number \(k\) is called the \emph{length} of the clause \(C_j\), and if for some \(k\) all clauses \(C_1,\dots,C_m\) have length at most \(k\), we say that \(C_1\land\dots\land C_m\) is a \(k\)-CNF representation. When only considering \(k\)-CNF representations for fixed \(k\), the problem is known as the \(k\)-SAT problem. For \(k=3\), \(n=7\), and \(m=4\) one could for instance consider the \(3\)-CNF formula
\[
f(x)=(x_1\lor\lnot x_3\lor x_4)\land (\lnot x_2\lor x_3)\land (\lnot x_1\lor x_6\lor\lnot x_7) \land x_7.
\]

To analyse the typical structure of the \(k\)-SAT problems we introduce the following probability distribution over the space of \(k\)-CNF formulas referred to as the  random \(k\)-SAT model \(F_k\).

\begin{definition}
For positive integers \(k,n,m\) we define \(F_k(n,m)\) to be the random \(k\)-CNF formula in \(n\) variables obtained as the conjunction of \(m\) i.i.d.\ sampled random clauses \(C_1,\dots,C_m\), each clause \(C_j\) being the disjunction of \(k\) elements chosen i.i.d.\ uniformly over the \(2n\) variables and their negations. 
\end{definition}

The random \(k\)-SAT model is often the one considered in the theoretical literature because of its simple structure. There is however more practical pertinence in analysing mixed models for random satisfiability; for instance, Gent and Walsh observed in \cite{GentWalsh1994} that satisfiability problems stemming from the industries most often are mixed (in the sense that the clauses are of varying lengths) and have a very different structure from that of \(F_k\). As a starting point we consider for distinct positive integers \(k_1\) and \(k_2\) the random mixed \(k_1\)- and \(k_2\)-SAT model \(F_{k_1,k_2}\). In this survey we focus on \(F_{1,2}\) and \(F_{2,3}\).

\begin{definition}
For all positive integers \(k_1,k_2,n,m_1,m_2\) we define the random CNF formula \(F_{k_1,k_2}(n,m_1,m_2)\) to be the conjunction of \(F_{k_1}(n,m_1)\) and \(F_{k_2}(n,m_2)\) defined above.
\end{definition}

\section{Conjectures}

\begin{figure}[h]
\centering
\includegraphics[scale=.67]{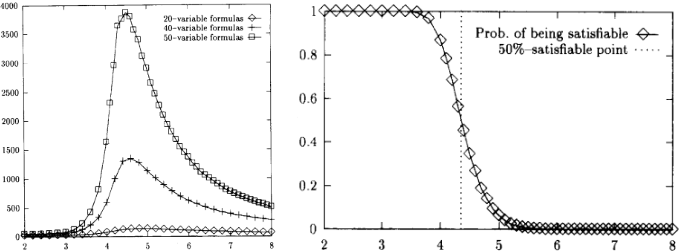}
\caption{Left: Fig.~2 from \cite{Selmanetal1996} showing median ``hardness'' as a function of \(\alpha\) for the model \(F_3(n,\alpha n)\) where \(n=20,40,50\). Right: Fig.~4 from \cite{Selmanetal1996} showing the empirical probability of satisfiability for \(F_3(50,50\alpha)\) as a function of \(\alpha\).}
\label{bosfig1}
\end{figure}

For almost a decade it was unclear how to produce formulas for which deciding satisfiability is hard. In '91, Cheeseman et al.\ publish the empirical study \cite{Cheesemanetal1991} where they produced hard satisfiability problems by transforming hard graph coloring problems (and the resulting CNF formulas are mixed). Soon after, Selman et al.\ found in \cite{Selmanetal1996} a way to directly produce hard instances using the random \(3\)-SAT model \(F_3(n,m)\). By considering the parameter \(\alpha=m/n\), the \emph{clause density}, they observed a spike in computational hardness in instances with clause density close to \(4.3\). This is shown on the left in Figure~\ref{bosfig1}. Next, they observed that the empirical probability of an \(F_3(n,\alpha n)\) instance being satisfiable, when \(\alpha\) varies from \(0\) to \(\infty\), drops rapidly from \(1\) to \(0\) at around the same point \(\alpha\approx 4.3\).

A sudden structural shift like the one suggested in Figure~\ref{bosfig1} when a parameter passes a single critical value is called a (sharp) \emph{phase transition}, and in 1992 it was famously conjectured by Chv{\'a}tal and Reed in their paper \cite{ChvatalReed1992} that the random \(k\)-SAT problem exhibits such a phenomenon (and proved in the same paper for the case \(k=2\); see current results in the next section):

\begin{conjecture}[\cite{ChvatalReed1992}]
\label{bosconj1}
For all \(k\geq 2\) there exists a value \(\alpha_\mathrm{sat}(k)>0\) separating the with high probability (w.h.p.) satisfiable instances of \(F_k(n,\alpha n)\) from the w.h.p.\ unsatisfiable ones, i.e.\
\[
\lim_{n\to\infty}\mathbb P(F_k(n,\alpha n)\text{ is satisfiable})=
\begin{cases}
1,&\text{if }\alpha<\alpha_\mathrm{sat}(k), \\
0,&\text{if }\alpha>\alpha_\mathrm{sat}(k).
\end{cases}
\]
\end{conjecture}

Selman et al.~\cite{Selmanetal1996} observed that there seems to be a close connection between the critical point of satisfiability \(\alpha_\mathrm{sat}(k)\) (if it exists) and the point where one finds hard instances of deciding satisfiability, and in the decades following commenced an intense research effort dedicated to proving Conjecture~\ref{bosconj1}, and this effort is still ongoing today. In what remains of this paper, we will  give an overview of the current status of progress towards Conjecture~\ref{bosconj1} and similar conjectures for random mixed models for satisfiability.

For the random mixed \(1\)- and \(2\)-SAT problem, no prior conjectures were put up (and the case \(k=1\) is suspiciously missing from Conjecture \ref{bosconj1}) apart from the following result shown in \cite{Bollobasetal2001} hinting at some larger theorem.

\begin{lemma}[\cite{Bollobasetal2001}]
\label{boslem1}
For all \(d>0\) it holds that
\[
\lim_{n\to\infty}\mathbb P(F_{1,2}(n,\log(n)^d,\alpha n)\text{ is satisfiable})=
\begin{cases}
1,&\text{if }\alpha<\alpha_\mathrm{sat}(2), \\
0,&\text{if }\alpha>\alpha_\mathrm{sat}(2).
\end{cases}
\]
\end{lemma}

The above result says that one can effectively ignore any polynomial in \(\log(n)\) \(1\)-clauses in the random \(2\)-SAT problem when \(n\to\infty\). An obvious question is then how many more \(1\)-clauses can be added before they become too numerous to be ignored.

For the \(F_{2,3}(n,\alpha_2 n,\alpha_3 n)\) model, it was conjectured in \cite{Monassonetal1996} that when fixing the proportion \(p\) of \(3\)-clauses in the mix, i.e.\ fixing \(p=\alpha_3/(\alpha_2+\alpha_3)\), then there again exists a critical point separating the w.h.p.\ satisfiable formulas from the w.h.p.\ unsatisfiable ones.

\begin{conjecture}[\cite{Monassonetal1996}]
\label{bosconj2}
For every \(p\in [0,1]\) there exists a value \(\alpha_\mathrm{sat}(2+p)>0\) such that for all \(\alpha_2,\alpha_3>0\) satisfying \(\alpha_3/(\alpha_2+\alpha_3)=p\) it holds that
\[
\lim_{n\to\infty}\mathbb P(F_{2,3}(n,\alpha_2 n,\alpha_3 n)\text{ is satisfiable})=
\begin{cases}
1,&\text{if }\alpha<\alpha_\mathrm{sat}(2+p), \\
0,&\text{if }\alpha>\alpha_\mathrm{sat}(2+p),
\end{cases}
\]
where \(\alpha=\alpha_2+\alpha_3\). Furthermore, for \(p<0.413\dots\) it holds that \(\alpha_\mathrm{sat}(2+p)=\alpha_\mathrm{sat}(2)/(1-p)\), and this does not hold for \(p\geq 0.413\dots\).
\end{conjecture}

Notice in the last part of  Conjecture~\ref{bosconj2} that \(\alpha<\alpha_\mathrm{sat}(2)/(1-p)\) if and only if \(\alpha_2<\alpha_\mathrm{sat}(2)\), which amounts to saying that if the proportion of \(3\)-clauses is low enough in \(F_{2,3}(n,\alpha_2 n,\alpha_3 n)\), then they can be ignored, and the problem behaves asymptotically as the \(F_2(n,\alpha_2 n)\) subformula. More concretely, one can add \(0.703 n\) random \(3\)-clauses to \(F_2(n,\alpha_2 n)\) for all \(\alpha_2<\alpha_\mathrm{sat}(2)\) and still have a w.h.p.\ satisfiable formula.

\section{Results}

Some serious progress has been made towards proving Conjecture~\ref{bosconj1}.

\begin{theorem}[\cite{ChvatalReed1992},\cite{Goerdt1996},\cite{Dingetal2022}]
\label{bosthm1}
Conjecture~\ref{bosconj1} holds for \(k=2\) with \(\alpha_\mathrm{sat}(2)=1\). Furthermore, there exists a positive integer \(k_0\) such that Conjecture~\ref{bosconj1} holds for all \(k\geq k_0\).
\end{theorem}

The first part of Theorem~\ref{bosthm1} was proved in the early '90s simultaneously by Chv\'{a}tal and Reed~\cite{ChvatalReed1992} and by Andreas Goerdt~\cite{Goerdt1996}. More is now known about the random \(2\)-SAT problem; see \cite{Bollobasetal2001}. The second part of the Theorem~\ref{bosthm1} was recently proved by Jian Ding, Allan Sly, and Nike Sun in \cite{Dingetal2022}. For the remaining \(k\) there are upper and lower bounds on \(\alpha_\mathrm{sat}(k)\) (if they exist) leaving only a gap of constant size (see \cite{Dingetal2022} and references therein).

In regards to Conjecture~\ref{bosconj2} we have the following.

\begin{theorem}[\cite{Achlioptasetal2001}]
\label{bosthm2}
Conjecture~\ref{bosconj2} holds for all \(p\leq 2/5\).
\end{theorem}

Achlioptas et al. \cite{Achlioptasetal2001} also give upper and lower bounds for \(\alpha_\mathrm{sat}(2+p)\) for the remaining \(p>2/5\).

Lastly, the random mixed \(1\)- and \(2\)-SAT problem has been solved by Andreas Basse-O'Connor, Tobias Overgaard, and Mette Skj{\o}tt in \cite{Connoretal2023}, thus ''completing'' Lemma~\ref{boslem1}.

\begin{theorem}[\cite{Connoretal2023}]
\label{bosthm3}
For all \(\alpha_1,q>0\) and \(0\leq\alpha_2<1\) it holds that
\[
\lim_{n\to\infty}\mathbb P(F_{1,2}(n,\alpha_1 n^q,\alpha_2 n)\text{ is satisfiable})=
\begin{cases}
1,&\text{if }q<1/2, \\
\exp\bigg(\frac{-\alpha_1^2}{4(1-\alpha_2)}\bigg),&\text{if }q=1/2, \\
0,&\text{if }q>1/2.
\end{cases}
\]
\end{theorem}

Theorem~\ref{bosthm3} states that one can add up to an order of \(\sqrt{n}\) random \(1\)-clauses to \(F_2(n,\alpha_2 n)\) for all \(\alpha_2<1\) and still have a w.h.p.\ satisfiable formula, and for the critical case \(q=1/2\) we even have an exact value for the limiting probability of satisfiability. Conditioning on the event that no \(1\)-clauses contradict each other only effects Theorem~\ref{bosthm3} in the case \(q=1/2\) where this limit is a different value. The case \(\alpha_2=0\) shows that Conjecture~\ref{bosconj1} does not hold for \(k=1\), firstly since one should consider in the order of \(\sqrt{n}\) \(1\)-clauses as opposed to order \(n\), and secondly since the random \(1\)-SAT problem does not even undergo a phase transition but instead has a smooth limiting probability of satisfiability.


\end{document}